\documentclass[11pt]{amsart}

\usepackage{amsmath,amssymb,amsfonts,amsthm,amscd,indentfirst,bigints,relsize}

\begin{document}

\def\hom{\mathop{\cH om\skp}}
\def\ext{\mathop{\cE xt\skp}}

\newtheorem{prop}{Proposition}[section]
\newtheorem{theo}{Theorem}[section]
\newtheorem{lemm}{Lemma}[section]
\newtheorem{coro}{Corollary}[section]
\newtheorem{rema}{Remark}[section]
\newtheorem{exam}{Example}[section]
\newtheorem{defi}{Definition}[section]
\newtheorem{conj}{Conjecture}[section]


\newtheorem{Def}{Definition}[section]
\newtheorem{Th}{Theorem}[section]
\newtheorem{Prop}{Proposition}[section]
\newtheorem{Not}{Notation}[section]
\newtheorem{Lemma}{Lemma}[section]
\newtheorem{Rem}{Remark}[section]
\newtheorem{Cor}{Corollary}[section]

\def\s{\section}                        \def\ss{\subsection}

\def\d{\begin{Def}}                     \def\t{\begin{Th}}
\def\p{\begin{Prop}}                    \def\n{\begin{Not}}
\def\la{\begin{Lemma}}                  \def\r{\begin{Rem}}
\def\c{\begin{Cor}}                     \def\ee{\begin{equation}}
\def\aa{\begin{eqnarray}}               \def\y{\begin{eqnarray*}}
\def\bd{\begin{description}}

\def\ed{\end{Def}}                      \def\et{\end{Th}}
\def\ep{\end{Prop}}                     \def\en{\end{No}}
\def\el{\end{Lemma}}                    \def\er{\end{Rem}}
\def\ec{\end{Cor}}                      \def\eee{\end{equation}}
\def\eaa{\end{eqnarray}}                \def\ey{\end{eqnarray*}}
\def\ebd{\end{description}}

\def\nn{\nonumber}                      \def\bp{{\bf Proof.}\hspace{2mm}}
\def\qe{\hfill{\rm Q.E.D.}}             \def\lj{\langle}
\def\rj{\rangle}                        \def\dd{\diamond}
\def\ox{\mbox{}}                        \def\lb{\label}

\def\begeq{\begin{equation}}
\def\endeq{\end{equation}}
\def\begarr{\begin{array}{rll}}
\def\endarr{\end{array}}
\def\earr{\begin{equation}\begin{array}{rll}}
\def\eearr{\end{array}\end{equation}}

\def\and{\quad{\rm and}\quad}
\def\dd{\bf $\diamond$}

\def\U{\mathcal{U}}
\def\V{\mathcal{V}}
\def\W{\mathcal{W}}
\def\R{\mathbb{R}}

\def\bl{\bigl(}
\def\br{\bigr)}
\def\lbe{_{\beta}}
\def\bN{{\mathbf N}}
\def\bs{{\mathbf s}}
\def\dist{{\mathbf dist}}
\def\<{\langle}
\def\>{\rangle}

\def\clo{\mathop{\rm cl\skp}}
\def\cdott{\!\cdot\!}
\def\cu{^{{\rm cu}}}

\def\lab{\label}

\def\Dint{\displaystyle\int}
\def\Dfrac{\displaystyle\frac}
\def\di{\displaystyle}

\title[Heat kernel Gaussian bounds on manifolds]{Heat kernel Gaussian bounds on manifolds I: manifolds with non-negative Ricci curvature}
 
\author{Xiangjin Xu }
\address{Department of Mathematical Sciences\\
         Binghamton University-SUNY\\
         Binghamton, New York, 13902, USA.}
\email{xxu@math.binghamton.edu}

\date{December 26th, 2019}

\begin{abstract}
This is first of series papers on new two-side Gaussian bounds for the heat kernel $H(x,y,t)$ on a complete manifold $(M,g)$. In this paper,  on a complete manifold $M$ with $Ric(M)\geq 0$, we obtain new two-side Gaussian bounds for the heat kernel $H(x,y,t)$, which improve the well-known Li-Yau's two-side bounds. As applications of our new two-side Gaussian bounds,  We obtain a sharp gradient estimate and a Laplacian estimate for the heat kernel on a complete manifold with $Ric(M)\geq 0$, and we also  give a simpler proof for the result concerning the asymptotic behavior in the time variable for the heat kernel as was proved in \cite{LiP-1} on a complete manifold $M$ with $Ric(M)\geq 0$ and maximal volume growth.
\end{abstract}

\maketitle

\tableofcontents

\section{\bf Introduction and main results}
Let $(M^n,g)$ be a complete Riemannian manifold with $Ric(M)\geq 0$. We say that $H(x, y, t)$,  defined on $M\times M\times(0, \infty)$, is a {\it heat kernel} if it is positive, symmetric in the x and y variables, and satisfies the heat equation
$$\left(\Delta_y-\frac{\partial}{\partial t}\right) H(x,y,t)=0\quad \text{with } \quad \lim_{t\rightarrow 0}H(x,y,t)=\delta_x(y),$$
where $\delta_x(y)$ denotes the point mass delta function at $x$ and $\Delta$ is the (negative definite) Laplace-Beltrami operator on $M$. For example, in $R^n$ the heat kernel is given by the classical Gauss-Weierstrass formula
\begin{eqnarray*}
H(x,y,t)=\left(\frac{1}{4\pi t}\right)^{n/2}e^{-\frac{d^2(x,y)}{4t}}.
\end{eqnarray*}
On a Riemannian manifold $(M,g)$, the local short time asymptotic expansion of the heat kernel $H(x,y,t)$ has been known for a long time (\cite{Ar},\cite{V}):\\
{\it There are smooth functions $H_i(x, y)$ defined on $(M \times M)\setminus C(M)$ with $C(M)=\{(x, y)| y \in Cut(x)\}$, 
\begin{equation}\label{Asymptotic}
H(x,y,t)\sim \Big(\frac{1}{4\pi t}\Big)^{n/2}e^{-\frac{d^2(x,y)}{4t}}\sum_{i=0}^{\infty}H_i(x,y) t^i
\end{equation}
holds uniformly as $t \rightarrow 0$ on any compact subsets of $(M \times M)\setminus C(M)$. 
}\\

A natural question is asked about the two-side bound estimates for the minimal heat kernel $H(x,y,t)$:\\
{\it {\bf Question:} Do there exist some nice functions $A(r,t)$ and $B(r,t)$, such that  the global two-side estimate holds
\begin{equation*}\label{TwoSideBound}
A(d(x,y),t)V_x^{-1}(\sqrt{t})e^{-\frac{d^2(x,y)}{4t}}\leq H(x,y,t)\leq B(d(x,y),t)V_x^{-1}(\sqrt{t}) e^{-\frac{d^2(x,y)}{4t}}
\end{equation*}
for all $(x,y,t)\in M\times M\times(0,\infty)$?
}\\

Dependence of the long time behavior of the heat kernel on the large scale geometry of $M$ is an interesting and important problem that has been intensively studied during the past few decades by many authors (see, for example, \cite{D}, \cite{G}, \cite{LTW}, \cite{LiP}, \cite{SY} and references therein). In the case of a complete manifold with $Ric(M)\geq 0$, in their pioneering work \cite{LY}, Li and Yau proved the following two side estimates, for all $x, y\in M$ and $t >0$,
\begin{equation}\label{LYBoundk=0}
C^{-1}(\delta)V_x^{-1}(\sqrt{t}) e^{-\frac{d^2(x,y)}{4(1-\delta)t}}\leq H(x,y,t)\leq C(\delta)V_x^{-1}(\sqrt{t}) e^{-\frac{d^2(x,y)}{4(1+\delta)t}}
\end{equation}
with $C(\delta)=c_1e^{\frac{c_2}{\delta}}$. Davies \cite{DE} proved an upper bound estimate of the heat hernel $(H(x,y,t)$ on a complete manifold $M$ with $Ric(M)\geq -k$ for  $r=\min\left\{1, \sqrt{t},\frac{t}{d}\right\}$ with $d=d(x,y)$:
\begin{equation}\label{DaviesBound}
0<H(x,y,t)\leq C(n,k)V_x^{-\frac{1}{2}}(r) V_y^{-\frac{1}{2}}(r)\exp\left[-\mu_1(M)t-\frac{d^2(x,y)}{4t}\right].
\end{equation} 
Unlike the Li-Yau\rq{}s estimates (\ref{LYBoundk=0}), Davies\rq{} estimate (\ref{DaviesBound}) only uses $V_x(r)$ with $r\leq 1$,  which provides some optimal upper bound estimates as $t\rightarrow 0$ as discussed at the end of paper \cite{DE}, but gives no decay as $t\rightarrow\infty$ even for on-diagonal estimates of $H(x,x,t)$. \\

In this paper, we follow the argument initially developed by Li and Yau in \cite{LY} with a new key observation that  one could replace $\delta t$ used in the proofs of the estimate (\ref{LYBoundk=0}) in \cite{LY} by some positive functions of $t$. Precisely, we show that:

\begin{Th}{\bf (Gaussian Lower Bounds)}\label{Thm-GLB-delta}
Let $(M^n,g)$ be a complete manifold with $Ric(M)\geq 0$, and $H(x, y, t)$ be its heat kernel, then for any $\delta>0$, we have the following Gaussian lower bound
\begin{eqnarray}\label{GLB-delta}
H(x,y,t)&\geq&e^{-\delta}\frac{ V_{\mathbb{R}^n}(R_{\delta}(t))}{V_x(R_{\delta}(t))}\left(4\pi t\right)^{-n/2}\exp\left[-\frac{d^2(x,y)}{4t}\right],
\end{eqnarray}
and the symmetrized version Gaussian lower bound
\begin{eqnarray}\label{GLB-delta-s}
H(x,y,t)&\geq&\frac{e^{-\delta} V_{\mathbb{R}^n}(R_{\delta}(t))}{\sqrt{V_x(R_{\delta}(t))V_y(R_{\delta}(t))}}\left(4\pi t\right)^{-n/2}\exp\left[-\frac{d^2(x,y)}{4t}\right],
\end{eqnarray}
with $d=d(x,y)$, and $R_{\delta}(t)=\frac{\sqrt{d^2+4\delta t}-d}{2}=\Big(\sqrt{\frac{d^2}{4\delta t}+1}+\sqrt{\frac{d^2}{4\delta t}}\Big)^{-1}\sqrt{\delta t}\leq \sqrt{\delta t}$.
\end{Th}
Furthermore, our lower bound estimates, Theorem \ref{Thm-GLB-delta}, are equivalent to the Cheeger and Yau\rq{}s heat kernel comparison theorem with $K=1$.

\begin{Th}{\bf (Gauss Upper Bound)}\label{Thm-GUB-delta}
Let $(M^n,g)$ be a complete manifold with $Ric(M)\geq 0$ and $H(x, y, t)$ be its heat kernel, then we have
\begin{eqnarray}\label{GUB-delta}
H(x,y,t)&\leq& f\left(\delta, \frac{d^2}{4t}\right)V_x^{-\frac{1}{2}}(R_{\delta}(t))V_y^{-\frac{1}{2}}(R_{\delta}(t))e^{-\frac{d^2(x,y)}{4t}}
\end{eqnarray}
and
\begin{eqnarray}\label{GUB-delta-S}
&&H(x,y,t)\leq e^{\delta}f^2\left(\delta, \frac{d^2}{4t}\right)\frac{(4\pi t)^{n/2}}{V_{\R^n}(R(t))} V_x^{-1}\left(R(t)\right)\exp \left[-\frac{d^2(x,y)}{4t}\right]
\end{eqnarray}
with $d=d(x,y)$ and
\begin{eqnarray}
&&R_{\delta}(t)=\frac{\sqrt{d^2+4\delta t}-d}{2}=\frac{2\delta t}{\sqrt{d^2+4\delta t}+d}
=\frac{\sqrt{\delta t}}{\sqrt{\frac{d^2}{4\delta t}+1}+\sqrt{\frac{d^2}{4\delta t}}},\nn\\
&&f\left(\delta, \frac{d^2}{4t}\right)=\left\{
\begin{array}{rl}
 e^{ \sqrt{\delta }+ \frac{\delta }{3}}\left(1+\sqrt{\delta}\right)^{\frac{n}{2}} , &  \frac{d^2}{4\delta t}\leq \frac{1}{3},\\
e^{2\delta} \left(1+\frac{4\delta t}{d^2}\right)^{\frac{n}{4}}, & \frac{d^2}{4\delta t}> \frac{1}{3}.
\end{array} \right.\nn
\end{eqnarray}
\end{Th}
If we choose $\delta=1$ in Theorem \ref{Thm-GLB-delta} and Theorem \ref{Thm-GUB-delta}, we could compare with the Li-Yau's two-side bound (\ref{LYBoundk=0}) with more details and discussion on Theorem \ref{GaussLBk=0Thm} and Theorem \ref{GaussUBk=0Thm}.\\

The first application of our two-side bound estimates is the following sharp gradient and  Laplacian estimates for the heat kernel $H(x,y,t)$ on a complete manifold with $Ric(M)\geq 0$,

\t\label{Hamilton-HK-sharp}
Suppose $(M^n,g)$ be a complete non-compact manifold with $Ric(M)\geq 0$, and $H(x,y,t)$ be its heat kernel, then  we have the following sharp gradient estimate for the heat kernel $H(x,y,t)$,
\begin{eqnarray}\label{HK-gradient-sharp}
&&\max\left\{\frac{1}{2},\frac{2}{1+\sqrt{1+\frac{4t}{d^2}}}\right\}t\left|\nabla\ln H(x,y,t)\right|^2\\
&\leq& C(n)+\frac{n}{2}\ln2+2n\ln\left(\sqrt{\frac{d^2}{4t}+1}+\sqrt{\frac{d^2}{4t}} \right)+\frac{d^2(x,y)}{4t},\nn
\end{eqnarray}
and the following Laplacian  estimate for the heat kernel $H(x,y,t)$,
\begin{eqnarray}\label{HK-Laplacian-sharp}
&&\max\left\{\frac{1}{2},\frac{2}{1+\sqrt{1+\frac{4t}{d^2}}}\right\}t\left(\frac{\Delta H(x,y,t)}{H(x,y,t)}\right)\\
&\leq& n+4C(n)+2n\ln2+8n\ln\left(\sqrt{\frac{d^2}{4t}+1}+\sqrt{\frac{d^2}{4t}} \right)+\frac{d^2(x,y)}{t},\nn
\end{eqnarray}
for all $x, y\in M$ and $t>0$, with 
$$C(n)=\frac{n}{2}\ln8\left(n+\sqrt{n^2+1}\right)+\ln \Gamma\left(\frac{n}{2}+1\right)+\frac{5-\sqrt{n^2+1}}{2}. $$
\et
\begin{Rem}
The sharpness of gradient estimate (\ref{HK-gradient-sharp}) achieves by the heat kernel of $\R^n$, $H(x,y,t)=\left(4\pi t\right)^{-n/2}e^{-\frac{d^2(x,y)}{4t}}$. As $\frac{d^2(x,y)}{4t}\rightarrow\infty$, the left-side of (\ref{HK-gradient-sharp}) asymptotically approaches to $\frac{d^2(x,y)}{4t}$, and the right-side  of (\ref{HK-gradient-sharp}) asymptotically approaches to $\frac{d^2(x,y)}{4t}+2n\ln\left(\sqrt{\frac{d^2}{4t}+1}+\sqrt{\frac{d^2}{4t}} \right)$. And the Laplacian estimate (\ref{HK-Laplacian-sharp}) is sharp in the order of $\frac{d^2(x,y)}{4t}$ for  the heat kernel of $\R^n$.\\
\end{Rem}  

Another application of our  two-side bound estimates is to give a simpler proof for the result concerning the asymptotic behavior in the time variable for the heat kernel as was proved in \cite{LiP} on the complete manifold $M$ with $Ric(M)\geq 0$ and maximal volume growth. We avoid the mess computations due to the $\delta$ loss in Li-Yau's two-side bounds (\ref{LYBoundk=0}) by our two-side exact Gaussian bounds, and we obtain the lower bound of the asymptotic behavior in the time variable for the heat kernel under much weaker condition, Lemma \ref{liminf-asym-Lemma}.\\
  
This paper is organized as following. In section 2, we list some Lemmas and Theorems from literature that we need to prove our main results. In section 3, we prove the two-side Gaussian bounds for the heat kernel $H(x,y,t)$ on a complete manifold with $Ric(M)\geq 0$. In section 4, as an application of our two-side bounds of the previous section, we yield estimates on the gradient and Laplacian of the heat kernel $H(x,y,t)$, which are sharp for the heat kernel on $\R^n$. In section 5, as another application, we give a simpler proof of  the asymptotic behavior in the time variable for the heat kernel \cite{LiP-1}.\\

In the forthcoming paper \cite{Xu}, we obtain two-side Gaussian bounds for the heat kernel $H(x,y,t)$ on a complete manifold $(M,g)$ with negative Ricci curvature lower bound, i.e. $Ric(M)\geq K$ for some $K>0$, using the improved Li-Yau's gradient estimates and induced parabolic Harnack inequalities
by J. F. Li and the author in \cite{LX}. Those two-side bound estimates improve the corresponding results by Li-Yau \cite{LY} and Sturm \cite{St}, by removing the $\delta$-loss in Gaussian term as we did in this paper on manifolds with non-negative Ricci curvature.\\ 

Throughout this paper, $(M,g)$ is assumed to be an n-dimensional complete connected Riemannian manifold with non-negative Ricci curvature, $Ric(M)\geq 0$, and denote $d=d_M(x,y)$ the shortest geodesic distance between $x$ and $y$ on $M$.\\

{\bf Acknowledgement:} Research of the author was supported in part by the NSF grant DMS-0852507, and Harpur College Grant in Support of Research, Scholarship and Creative Work in Year 2010-11, 2012-13, and Harpur Faculty Research Award in Year 2019-20.\\

\section{\bf Some preliminaries}

In this section, we list some Lemmas and Theorems from literature that we need to prove our main results. The two side bound estimates (\ref{LYBoundk=0}) are based on the following Li-Yau\rq{}s gradient estimates and induced parabolic Harnack inequalities proved in Li and Yau\rq{}s pioneering work \cite{LY}:

\t{\bf (Li-Yau \cite{LY}).} \label{LYThm}
Let $(M,g)$ be a complete Riemannian manifold with $Ric(M)\geq 0$. If $u(x,t)$ is a positive solution of the heat equation $(\partial_t-\Delta)u(x,t)=0$, there is the sharp gradient estimate: 
\begin{eqnarray}\label{LYGradient}
\frac{\left|\nabla u\right|^2}{u^2}-\frac{u_t}{u}\leq \frac{ n}{2t},
\end{eqnarray}
and the Harnack inequality:
\begin{eqnarray}\label{LYHarnack}
\frac{u(x,t_1)}{u(y,t_2)}\leq \left(\frac{t_2}{t_1}\right)^{\frac{ n}{2}}\exp\left[\frac{ d^2(x,y)}{4(t_2-t_1)}\right],
\end{eqnarray}
for any $x,y\in M$ and $0<t_1<t_2<\infty$.
\et

\c\label{MVIneq}{\bf (Mean Value Inequalities)}
Under the same assumption as in Theorem \ref{LYThm} above, then the following mean value inequality holds: for $0<R<\infty$ and $0<t_1<t_2<\infty$, 
\begeq\label{MVIneq1} 
u(x,t_1)\le \left(\oint_{B_{x}(R)} u^p(y,t_2)dy\right)^{\frac{1}{p}}\left(\frac{t_2}{t_1}\right)^{n/2}\cdot \exp\left[\frac{R^2}{4(t_2-t_1)}\right]
\endeq
where $p>0$, and $\oint_{B_{x}(R)}=Vol^{-1}(B_x(R))\cdot\int_{B_{x}(R)}.$
\ec  

And the Cheeger and Yau\rq{}s heat kernel comparison theorem:

\t{\bf (Cheeger-Yau \cite{CY})}\label{Cheeger-Yau}
Let $M$ be a complete manifold without boundary with $Ric(M)\geq -(n-1)K$ and $H(x, y, t)$ denotes the minimal heat kernel defined on $M$ and $H^K(x,y,t)$ be the heat kernel on the simply connected space form $M^K$ with constant sectional curvature $-K \leq 0$. Then
\begin{eqnarray}\label{Cheeger-Yau-global}
H(x,y,t)\geq H^K(\bar{x},\bar{y}, t)
\end{eqnarray}
are valid for $d_M(x,y)=d_{M^K}(\bar{x},\bar{y})$ and for $t\in (0,\infty)$.
\et

We need the following Lemma due to Davies \cite{D} to estimate the integral of the heat kernel $H(x, y, t)$,

\la\label{Davies} {\bf (Lemma 13.2 (Davies) in \cite{LiP})}
Let $(M^n,g)$ be a complete manifold with $Ric(M)\geq 0$.  Suppose $B_1$ and $B_2$ are bounded subsets of $M$,
$$\int_{B_1}\int_{B_2} H(x, y, t) dy dx \leq V^{\frac{1}{2}}(B_1) V^{\frac{1}{2}}(B_2) \exp\left(-\frac{d^2(B_1,B_2)}{4t}\right),$$
where $V (B_i)$ denotes the volume of the set $B_i$ for $i = 1, 2$ and $d(B_1, B_2)$ denotes the distance between the sets $B_1$ and $B_2$. 
\el

We'll use the following Hamilton's gradient estimates,  closed manifolds case by Hamilton \cite{H} and complete manifolds case by  Kotschwar \cite{K}:

\t\label{H-gradient} {\bf (Hamilton \cite{H}, Kotschwar \cite{K}).}
Suppose $(M^n,g)$ be a complete manifold with $Ric(M)\geq -K$ for some constant $K\geq 0$. If $0<u(x,t)\leq A$ for some constant $0<A<\infty$  is a bounded positive solution to the heat equation on $M\times[0,T]$ for $0<T\leq \infty$, then
$$t\left|\nabla\ln u(x,t)\right|^2\leq (1+2Kt)\ln\left(\frac{A}{u(x,t)}\right),\quad \forall\; (x,t)\in M\times[0,T].$$
\et

And the following Hamilton's Laplacian estimates, closed manifolds case by Hamilton \cite{H} and   complete manifolds case by J.Y. Wu \cite{W}:

\t\label{H-Laplacian} {\bf (Hamilton \cite{H}, J. Y. Wu \cite{W}).}
Suppose $(M^n,g)$ be a complete manifold with $Ric(M)\geq 0$. If $0<u(x,t)\leq  A$ for some constant $0<A<\infty$ is a bounded positive solution to the heat equation on $M\times[0,T]$ for $0<T\leq \infty$, then
$$t\left(\frac{\Delta u(x,t)}{u(x,t)}\right)\leq n+4\ln\left(\frac{A}{u(x,t)}\right),\quad \forall\; (x,t)\in M\times[0,T].$$
\et

\section{\bf Two side Gaussian bound estimates of $H(x,y,t)$}

In this section, we follow the main  argument in Chapter 13 in Li \cite{LiP}, much of which was first developed by Li and Yau \cite{LY}. A new key observation is that one could replace $\delta t$ used in Theorem 13.1 and 13.3 in \cite{LiP}, by some positive continuous functions of $t$.

\subsection{\bf Lower bound estimates of $H(x,y,t)$} 
Based on the Harnack inequality (\ref{LYHarnack}),  we will derive a general lower bound for the heat kernel $H(x,y,t)$:

\p\label{Prop-LB}
Let $(M^n,g)$ be a complete manifold with  $Ric(M)\geq 0$, $H(x, y, t)$ be its heat kernel, then for all $0<T<t$, we have 

\begin{eqnarray}\label{LBpre}
H(x,y,t)&\geq& \left(\frac{T}{t}\right)^{n/2}\exp\left[-\frac{R^2}{t-T}\right]\frac{V_{\R^n}(R)}{V_x(R)} \left(\frac{1}{4\pi T}\right)^{n/2}e^{-\frac{d^2(x,y)}{4T}},
\end{eqnarray}
where $0<T=T(t)<t$ is some positive continuous function on $t\in (0,\infty)$ and $R=R(t)$ is any positive continuous increasing function on $t\in (0,\infty)$. 
\ep

\bp
For a fixed $y\in M$, applying Corollary \ref{MVIneq} on $B_x(R)$ to the solution $f(x, t) = H(x, y, t)$ by taking  $0<t_0<t$ and $z\in B_x(R)$, we have
\begin{eqnarray*}
\int_{B_x(R)}H(z,y,t_0) dz\leq H(x,y,t)V_x(R) \left(\frac{t}{t_0}\right)^{n/2}\exp\left[\frac{R^2}{4(t-t_0)}\right].
\end{eqnarray*} 
Applying Corollary \ref{MVIneq} on $B_{\bar{x}}(R)$ to the solution $f(\bar{x}, s) = \left(4\pi s\right)^{-\frac{n}{2}}e^{-\frac{d^2(\bar{x},\bar{y})}{4s}}$, the heat kernel on $\R^n$ with  $d_M(x,y)=d_{\R^n}(\bar{x},\bar{y})$, by taking $0<T< t_0$, we obtain

\begin{eqnarray*}
f(\bar{x}, T)V_{\bar{x}}(R)\leq \left(\frac{t_0}{T}\right)^{n/2}\exp\left[\frac{R^2}{4(t_0-T)}\right]\int_{B_{\bar{x}}(R)}f(\bar{z}, t_0)d\bar{z}.
\end{eqnarray*} 
Using the Cheeger-Yau's Heat Kernel comparison theorem \ref{Cheeger-Yau} with $K=0$, 
$$ H(z,y,t_0)\geq \left(\frac{1}{4\pi t_0}\right)^{n/2}e^{-\frac{d^2(z,y)}{4t_0}}=f(\bar{z},t_0)$$
are valid for $d_M(z,y)=d_{\R^n}(\bar{z},\bar{y})$ and for $t\in (0,\infty)$. Let $t-t_0=t_0-T$, we obtain (\ref{LBpre}) from above estimates.
\qe

By choosing suitable functions $R(t)>0$ and $T(t)<t$ in Proposition \ref{Prop-LB}, we will prove Theorem \ref{Thm-GLB-delta} as following: \\

{\bf Proof of Theorem \ref{Thm-GLB-delta}:}
Firstly for $d=d(x,y)>0$,  from (\ref{LBpre}), we obtain
\begin{eqnarray*}
&&H(x,y,t)\\
&\geq& \left(\frac{T}{t}\right)^{\frac{n}{2}}\exp\left[-\frac{R^2}{t-T}\right] V_x^{-1}(R)V_{\mathbb{R}^n}(R) \left(4\pi T\right)^{-\frac{n}{2}}\exp\left[-\frac{d^2(x,y)}{4T}\right]\\
&=& V_{\mathbb{R}^n}(R)V_x^{-1}(R)\exp\left[-\frac{R^2}{t-T}-\frac{d^2(t-T)}{4tT}\right]\left(4\pi t\right)^{-\frac{n}{2}}\exp\left[-\frac{d^2(x,y)}{4t}\right].
\end{eqnarray*}
Let $T(t)=\frac{d}{d+2R}t<t$, then we have
\begin{eqnarray*}
H(x,y,t)&\geq& V_{\mathbb{R}^n}(R)V_x^{-1}(R)\exp\left[-\frac{dR+R^2}{t}\right]\left(4\pi t\right)^{-\frac{n}{2}}\exp\left[-\frac{d^2(x,y)}{4t}\right]
\end{eqnarray*}
Next to choose $\frac{dR+R^2}{t}=\delta>0$, which gives 
\begin{eqnarray}\label{RandT}
\left\{
\begin{array}{rl}
R_{\delta}(t)&=\frac{\sqrt{d^2+4\delta t}-d}{2}=\frac{2\delta t}{\sqrt{d^2+4\delta t}+d}
=\frac{\delta \sqrt{t}}{\sqrt{\frac{d^2}{4t}+\delta}+\sqrt{\frac{d^2}{4t}}}, \\
T_{\delta}(t)&=\frac{d}{d+ 2 R_{\delta}(t)}t=\frac{d}{\sqrt{d^2+4\delta t}}t<t,\quad d>0.
\end{array} \right.
\end{eqnarray}
We obtain
\begin{eqnarray*}
H(x,y,t)&\geq&e^{-\delta}\frac{ V_{\mathbb{R}^n}(R_{\delta}(t))}{V_x(R_{\delta}(t))}\left(4\pi t\right)^{-\frac{n}{2}}\exp\left[-\frac{d^2(x,y)}{4t}\right].
\end{eqnarray*}
The above estimate still holds for $d=d_M(x,y)=0$ by the continuity of $H(x,y,t)$, which proves the estimate (\ref{GLB-delta}).

\qe

Next we show that Theorem \ref{Thm-GLB-delta} will imply the Cheeger-Yau's heat kernel comparison Theorem with $K=0$. More precisely, we have

\t\label{LB->CY} 
Under the same assumptions of Theorem \ref{Thm-GLB-delta},  the lower bound estimates (\ref{GLB-delta}) with any $\delta>0$ imply the estimate (\ref{Cheeger-Yau-global}) with $K=0$ in Theorem \ref{Cheeger-Yau}.
\et

\bp
For fixed $t>0$ and $d=d_M(x,y)>0$, there are
$$\lim_{\delta\rightarrow 0}R_{\delta}(t)=0,\quad \quad \lim_{\delta\rightarrow 0}\frac{V_{\R^n}(R_{\delta}(t))}{V_x(R_{\delta}(t))}=\lim_{r\rightarrow 0}\frac{V_{\R^n}(r)}{V_x(r)}=1.  $$
Hence let $\delta\rightarrow 0$ in the estimates (\ref{GLB-delta}), one has
\begin{eqnarray}\label{d>0}
H(x,y,t)\geq \left(\frac{1}{4\pi t}\right)^{n/2}e^{-\frac{d^2(x,y)}{4t}},
\end{eqnarray}
for any given $t>0$ and $x\neq y$, which implies 
\begin{eqnarray}\label{d=0}
H(x,x,t)\geq \left(\frac{1}{4\pi t}\right)^{n/2},
\end{eqnarray}
by the continuity of the heat kernels when $t>0$. Combining (\ref{d>0}) and (\ref{d=0}), the estimate (\ref{Cheeger-Yau-global}) with $K=0$ in  Theorem \ref{Cheeger-Yau} follows.
\qe\\

Setting $\delta=1$ in (\ref{GLB-delta}), we have the following Gaussian lower bounds of the heat kernel $H(x,y,t)$ with precise constants,

\t{\bf (Gaussian Lower Bounds)}\label{GaussLBk=0Thm}
Let $(M^n,g)$ be a complete manifold with $Ric(M)\geq 0$, and $H(x, y, t)$be its heat kernel, then we have the following Gaussian lower bound
\begin{eqnarray}\label{GLBk=0-1}
&&H(x,y,t)
\geq \frac{ V_{\mathbb{R}^n}(R(t))}{eV_x(R(t))}\left(4\pi t\right)^{-\frac{n}{2}}\exp\left[-\frac{d^2(x,y)}{4t}\right]\\
&\geq& \frac{ V_{\mathbb{R}^n}(1)}{e(4\pi)^{n/2}}\left(\sqrt{\frac{d^2}{4t}+1}+\sqrt{\frac{d^2}{4t}} \right)^{-n}V_x^{-1}\left( \sqrt{t}\right)\exp\left[-\frac{d^2(x,y)}{4t}\right],\nn
\end{eqnarray}
and  the symmetrized version 
\begin{eqnarray}\label{GLBk=0-2}
&&\quad H(x,y,t)
\geq \frac{e^{-1} V_{\mathbb{R}^n}(R(t))}{\sqrt{V_x(R(t))V_y(R(t))}}\left(4\pi t\right)^{-\frac{n}{2}}\exp\left[-\frac{d^2(x,y)}{4t}\right]\\
&\geq& \frac{ V_{\mathbb{R}^n}(1)}{e(4\pi)^{n/2}}\left(\sqrt{\frac{d^2}{4t}+1}+\sqrt{\frac{d^2}{4t}} \right)^{-n}V_x^{-\frac{1}{2}}\left(\sqrt{t}\right)V_y^{-\frac{1}{2}}\left(\sqrt{t}\right)\exp\left[-\frac{d^2(x,y)}{4t}\right].\nn
\end{eqnarray}
with $d=d(x,y)$, and $R(t)=\frac{\sqrt{d^2+4t}-d}{2}=\Big(\sqrt{\frac{d^2}{4t}+1}+\sqrt{\frac{d^2}{4t}}\Big)^{-1}\sqrt{t}\leq \sqrt{t}$. \\
\et

\begin{Rem}
Comparing with Li-Yau\rq{}s lower bound estimate (\ref{LYBoundk=0}), after optimizing with $\delta>0$, which have an exponential lost, $exp\left[-C\sqrt{\frac{d^2}{t}}\right]$, to the exact  Gaussian bound $V_x^{-1}(\sqrt{t})\exp\left[-\frac{d^2(x,y)}{4t}\right]$. But our lower bound estimate (\ref{GLBk=0-1}) has only a polynomial lost, $\left(\sqrt{\frac{d^2}{4t}+1}+\sqrt{\frac{d^2}{4t}}\right)^{-n}$, to the exact Gaussian bound $V_x^{-1}(\sqrt{t})\exp\left[-\frac{d^2(x,y)}{4t}\right]$.
\end{Rem}

\begin{Rem}\label{Rem on GLBk=0-2}
When  the manifold $M$, such as with end as $R^{\tau}\times S^{n-\tau}$, has slower volume growth as $\lim_{R\rightarrow\infty}\frac{V_x(R)}{R^{\tau}}= C_M>0$ for some $0\leq \tau<n$, estimates (\ref{GLBk=0-1}) and (\ref{GLBk=0-2}) provide,
\begin{eqnarray*}
H(x,y,t)\gtrsim \frac{V_{\mathbb{R}^n}(1)}{e C_M}R(t)^{n-\tau}\big(4\pi t\big)^{-n/2}\exp\left[-\frac{d^2(x,y)}{4t}\right], \text{ as } t\rightarrow\infty,
\end{eqnarray*}
which will provide  better lower bounds than the Gauss kernel of $R^n$.
\end{Rem}

\subsection{\bf Upper bound estimates of $H(x,y,t)$} 
Follow the argument first developed by Li and Yau \cite{LY}, we first give the proof of  Theorem \ref{Thm-GUB-delta}:\\

{\bf Proof of Theorem \ref{Thm-GUB-delta}:}
Fixed $y\in M$, applying Corollary \ref{MVIneq} on $B_x(R)$ to the solution $f(x, t) = H(x, y, t)$ by taking $t_1 = t$ and $t_2 = T_0>t$, we have
\begin{eqnarray*}\label{T0(t)}
H(x,y,t)\leq \left(\frac{T_0}{t}\right)^{n/2}\exp\left[\frac{R^2}{4(T_0-t)}\right]V_x^{-1}(R)\int_{B_x(R)} H(x',y,T_0)dx'.
\end{eqnarray*} 
Applying Corollary \ref{MVIneq} on $B_y(R)$ to the solution $f(x, s) = \int_{B_x(R)}H(x', y, s)dx'$ by taking $t_1 = T_0$ and $t_2 = T>T_0$, we obtain
\begin{eqnarray*}\label{T(t)}
\int_{B_x(R)}H(x',y,T_0)dx'&\leq& \left(\frac{T}{T_0}\right)^{n/2}\exp\left[\frac{R^2}{4(T-T_0)}\right] V_y^{-1}(R)\\
&&\times \int_{B_y(R)}\int_{B_x(R)} H(x',y',T)dx'dy'.
\end{eqnarray*} 
Hence combining the above estimates, let $T(t)-T_0(t)=T_0(t)-t$, we have
\begin{eqnarray*}\label{H(x,y,t)pre}
H(x,y,t)&\leq& \left(\frac{T}{t}\right)^{n/2}\exp\left[\frac{R^2}{T-t}\right]V_x^{-1}(R)V_y^{-1}(R)\\
&&\times \int_{B_y(R)}\int_{B_x(R)} H(x',y',T)dx'dy'\nn
\end{eqnarray*}
On the other hand, Lemma \ref{Davies} implies that
\begin{eqnarray*}\label{DaviesIneq}
\int_{B_y(R)}\int_{B_x(R)} H(x',y',T)dx'dy'\leq V_x^{\frac{1}{2}}(R)V_y^{\frac{1}{2}}(R) \exp \left[-\frac{d^2(B_x(R),B_y(R))}{4T}\right].
\end{eqnarray*}
Combing the above two estimates, we obtain
\begin{eqnarray*}
H(x,y,t)\leq \left(\frac{T}{t}\right)^{n/2}V_x^{-\frac{1}{2}}(R)V_y^{-\frac{1}{2}}(R)\exp\left[\frac{R^2}{T(t)-t}-\frac{d^2(B_x(R),B_y(R))}{4T(t)}\right].
\end{eqnarray*}
Observing that
\begin{equation*}
d(B_x(R),B_y(R))=\left(d(x,y)-2R\right)_+=\left\{\begin{array}{ccc}
0, & if\; d(x,y)\leq 2R\\
d(x,y)-2R, & if\; d(x,y)>2R,
\end{array} \right .
\end{equation*}
we have
\begin{eqnarray}\label{Ondiagonal}
H(x,y,t)&\leq& \left(\frac{ T(t)}{t}\right)^{\frac{n}{2}} \exp\left[\frac{R^2(t)}{T(t)-t}+\frac{d^2}{4t}-\frac{\left(d-2R(t)\right)_+^2}{4T(t)}\right]\\
&&\times V_x^{-\frac{1}{2}}(R(t))V_y^{-\frac{1}{2}}(R(t))exp\left[-\frac{d^2(x,y)}{4t}\right],\nn
\end{eqnarray}
where $R(t)>0$ and $T(t)>t$ are the positive functions to be chosen as the following: for each $\delta>0$, define 
\begin{equation}\label{RandT-UB}
\left\{
\begin{array}{rl}
R_{\delta}(t)&=\frac{\sqrt{d^2+4\delta t}-d}{2}=\frac{2\delta t}{\sqrt{d^2+4\delta t}+d}
=\frac{\delta\sqrt{t}}{\sqrt{\frac{d^2}{4t}+\delta}+\sqrt{\frac{d^2}{4t}}}, \\
&\\
T_{\delta}(t)&=\left\{
\begin{array}{rl}
\left(1+\sqrt{\delta}\right)t>t, & \text{if} \quad \frac{d^2}{4\delta t}\leq \frac{1}{3},\\
\left(\sqrt{1+\frac{4\delta t}{d^2}}\right)t>t, & \text{if} \quad \frac{d^2}{4\delta t}> \frac{1}{3},
\end{array} \right.
\end{array} \right.
\end{equation}

{\bf Case 1:} when $0\leq d=d(x,y)\leq 2R_{\delta}(t)$, equivalent to $0\leq \frac{d^2}{4t}\leq \frac{\delta}{3}$, we have $T(t)=\left(1+\sqrt{\delta}\right)t>t$, and
$$\frac{R_{\delta}^2}{T_{\delta}-t}+\frac{d^2}{4t}-\frac{\left(d-2R_{\delta}\right)_+^2}{4T_{\delta}}=\left(\sqrt{1+\frac{d^2}{4\delta t}}-\sqrt{\frac{d^2}{4\delta t}}\right)^2\sqrt{\delta }+\frac{d^2}{4t}\leq \sqrt{\delta }+ \frac{\delta }{3}. $$
Hence we have
\begin{equation}\label{GUB-I}
H(x,y,t)\leq \left(1+\sqrt{\delta}\right)^{\frac{n}{2}} e^{ \sqrt{\delta }+ \frac{\delta }{3}} V_x^{-\frac{1}{2}}(R_{\delta}(t))V_y^{-\frac{1}{2}}(R_{\delta}(t))e^{-\frac{d^2(x,y)}{4t}}.
\end{equation}

{\bf Case II:}  when $ d=d(x,y)> 2R_{\delta}(t)$, equivalent to $\frac{d^2}{4t}> \frac{\delta}{3}$, we have $T(t)=\left(\sqrt{1+\frac{4\delta t}{d^2}}\right)t>t$, and
\begin{eqnarray*}
\frac{R_{\delta}^2}{T_{\delta}-t}+\frac{d^2}{4t}-\frac{\left(d-2R_{\delta}\right)_+^2}{4T_{\delta}}&=&\frac{R_{\delta}^2}{T_{\delta}-t}+\frac{d^2}{4t}-\frac{d^2-4dR_{\delta}+4R_{\delta}^2}{4T_{\delta}}\\
&=&\frac{2R_{\delta} d}{T_{\delta}}=\frac{4\delta}{\left(1+\sqrt{1+\frac{4\delta t}{d^2}}\right)\sqrt{1+\frac{4\delta t}{d^2}}}\leq 2\delta . 
\end{eqnarray*}\\
Hence we have
\begin{eqnarray}\label{GUB-II}
H(x,y,t)&\leq& e^{2\delta} \left(1+\frac{4\delta t}{d^2}\right)^{\frac{n}{4}} V_x^{-\frac{1}{2}}(R_{\delta}(t))V_y^{-\frac{1}{2}}(R_{\delta}(t))e^{-\frac{d^2(x,y)}{4t}}.
\end{eqnarray}

Combining {\bf Case I} and {\bf Case II}, we have the estimate (\ref{GUB-delta}) as
\begin{eqnarray*}
H(x,y,t)\leq f\left(\delta, \frac{d^2}{4t}\right)  V_x^{-\frac{1}{2}}(R(t))V_y^{-\frac{1}{2}}(R(t))exp\left[-\frac{d^2(x,y)}{4t}\right].
\end{eqnarray*}
On the other hand, from  (\ref{GLB-delta}) we have
\begin{eqnarray*}
H(x,y,t)&\geq& e^{-\delta}\frac{ V_{\mathbb{R}^n}(R_{\delta}(t))}{V_y(R_{\delta}(t))}\left(4\pi t\right)^{-n/2}\exp\left[-\frac{d^2(x,y)}{4t}\right],
\end{eqnarray*}
therefore,
\begin{eqnarray*}
V_y^{-\frac{1}{2}}\left(R_{\delta}(t)\right)\leq  e^{\delta}f\left(\delta, \frac{d^2}{4t}\right)\frac{\left(4\pi t\right)^{n/2}}{ V_{\R^n}(R_{\delta}(t))}V_x^{-\frac{1}{2}}\left(R_{\delta}(t)\right),
\end{eqnarray*}
substituting this in (\ref{GUB-delta}), we get
\begin{eqnarray*}
&&H(x,y,t)\leq e^{\delta}f^2\left(\delta, \frac{d^2}{4t}\right)\frac{\left(4\pi t\right)^{n/2}}{ V_{\R^n}(R_{\delta}(t))}V_x^{-1}\left(R_{\delta}(t)\right)\exp \left[-\frac{d^2(x,y)}{4t}\right].
\end{eqnarray*}

\qe

By choosing $\delta=1$ in the  estimate (\ref{GUB-delta}) and and applying Bishop's volume comparison theorem \cite{BC}, we have the following Gaussian upper bound for the heat kernel $H(x,y,t)$ with precise constants,

\t{\bf (Gauss Upper Bound)}\label{GaussUBk=0Thm}
Let $(M^n,g)$ be a complete manifold with $Ric(M)\geq 0$, and $H(x, y, t)$ be its heat kernel, then  we have
\begin{eqnarray}\label{GUBk=0}
&&H(x,y,t)
\leq 2^{\frac{n}{2}} e^{2} V_x^{-\frac{1}{2}}(R(t))V_y^{-\frac{1}{2}}(R(t))exp\left[-\frac{d^2(x,y)}{4t}\right]\\
&\leq& 2^{\frac{n}{2}} e^{2}\left(\frac{\sqrt{t}}{R(t)}\right)^{n}V_x^{-\frac{1}{2}}\left(\sqrt{t}\right)V_y^{-\frac{1}{2}}\left(\sqrt{t}\right)\exp \left[-\frac{d^2(x,y)}{4t}\right]\nn\\
&=&2^{\frac{n}{2}} e^{2}\left(\sqrt{\frac{d^2}{4t}+1}+\sqrt{\frac{d^2}{4t}}\right)^{n}V_x^{-\frac{1}{2}}\left(\sqrt{t}\right)V_y^{-\frac{1}{2}}\left(\sqrt{t}\right)\exp \left[-\frac{d^2(x,y)}{4t}\right],\nn
\end{eqnarray}
and
\begin{eqnarray}\label{GUBk=0S}
&&H(x,y,t)\leq \frac{e^5 2^n(4\pi t)^{n/2}}{V_{\R^n}(R(t))} V_x^{-1}\left(R(t)\right)\exp \left[-\frac{d^2(x,y)}{4t}\right]\\
&\leq& \frac{e^5 (16\pi)^{n/2}}{V_{\R^n}(1)}\left(\sqrt{\frac{d^2}{4t}+1}+\sqrt{\frac{d^2}{4t}}\right)^{n} V_x^{-1}\left(R(t)\right)\exp \left[-\frac{d^2(x,y)}{4t}\right]\nn\\
&\leq& \frac{e^5 (16\pi)^{n/2}}{V_{\R^n}(1)}\left(\sqrt{\frac{d^2}{4t}+1}+\sqrt{\frac{d^2}{4t}}\right)^{2n} V_x^{-1}\left(\sqrt{t}\right)\exp \left[-\frac{d^2(x,y)}{4t}\right].\nn
\end{eqnarray}
with $d=d(x,y)$, and $R(t)=\frac{\sqrt{d^2+4t}-d}{2}=\Big(\sqrt{\frac{d^2}{4t}+1}+\sqrt{\frac{d^2}{4t}}\Big)^{-1}\sqrt{t}\leq \sqrt{t}$. 
\et

\begin{Rem}
Comparing with Li-Yau\rq{}s upper bound estimate (\ref{LYBoundk=0}), after optimizing with $\delta>0$, which have an exponential lost, $exp\left[C\sqrt{\frac{d^2}{t}}\right]$, to the exact  Gaussian bound $V_x^{-1/2}(\sqrt{t}) V_y^{-1/2}(\sqrt{t})\exp\left[-\frac{d^2(x,y)}{4t}\right]$. But our upper bound estimate (\ref{GUBk=0}) have only a polynomial lost, $\left(\sqrt{\frac{d^2}{4t}+1}+\sqrt{\frac{d^2}{4t}}\right)^{n}$, to the exact Gaussian bound $V_x^{-1/2}(\sqrt{t}) V_y^{-1/2}(\sqrt{t})\exp\left[-\frac{d^2(x,y)}{4t}\right]$.
\end{Rem}

\begin{Rem}
When the manifold $M$, such as with end as $R^{\tau}\times S^{n-\tau}$, has slower volume growth as $\lim_{R\rightarrow\infty}\frac{V_x(R)}{R^{\tau}}= C_M>0$ for some $0\leq \tau<n$, the  estimate (\ref{GUBk=0}) provide, as $t\rightarrow\infty$
\begin{eqnarray*}
H(x,y,t)&\lesssim& C\left(R(t)\right)^{-\tau}\exp\left[-\frac{d^2(x,y)}{4t}\right]\\
&=& C\left(\sqrt{\frac{d^2}{4t}+1}+\sqrt{\frac{d^2}{4t}}\right)^{\tau}t^{-\tau/2}\exp\left[-\frac{d^2(x,y)}{4t}\right] ,
\end{eqnarray*}
which is better than the Gaussian heat kernel of $\R^n$.\\
\end{Rem}

\section{\bf Gradient and Laplacian estimates for $H(x,y,t)$}

In this section, an application of our two side bound estimates of the heat kernel $H(x,y,t)$ of the previous section, together with the gradient estimate (Theorem \ref{H-gradient}) and Laplacian estimate (Theorem \ref{H-Laplacian}), yields estimates on the gradient and Laplacian of the heat kernel $H(x,y,t)$ for complete manifolds with $Ric(M)\geq 0$ that is sharp for the heat kernel on $\R^n$. To do these, we followthe arguments in \cite{K} and in \cite{W}, where they made use of Li-Yau's two-side bound (\ref{LYBoundk=0}) \cite{LY}, which has a $\delta$-loss in Gaussian term and the constant $C(\delta)$ blows up as $\delta\rightarrow 0$. Precisely, we firstly show

\t\label{Hamilton-HK}
Suppose $(M^n,g)$ be a complete non-compact manifold with $Ric(M)\geq 0$, and $H(x,y,t)$ be its heat kernel, then for any $0<\alpha<1$, we have the following gradient estimate
\begin{eqnarray*}
&&\left(1-\alpha\right)t\left|\nabla\ln H(x,y,t)\right|^2\\
&\leq& C(n)-\frac{n}{2}\ln\alpha+n\ln\left(\sqrt{\frac{d^2}{4t}+1}+\sqrt{\frac{d^2}{4t}} \right)+\frac{d^2(x,y)}{4t},
\end{eqnarray*}
and the following Laplacian estimate
\begin{eqnarray*}
&&\left(1-\alpha\right)t\left(\frac{\Delta H(x,y,t)}{H(x,y,t)}\right)\\
&\leq& n+4C(n)-2n\ln\alpha+4n\ln\left(\sqrt{\frac{d^2}{4t}+1}+\sqrt{\frac{d^2}{4t}} \right)+\frac{d^2(x,y)}{t},
\end{eqnarray*}
for all $x, y\in M$ and $t>0$, with 
$$C(n)=\frac{n}{2}\ln8\left(n+\sqrt{n^2+1}\right)+\ln \Gamma\left(\frac{n}{2}+1\right)+\frac{5-\sqrt{n^2+1}}{2}. $$
\et

\bp
From our Gaussian lower bound (\ref{GLBk=0-2}), we have
\begin{eqnarray}\label{HK-LB}
\quad H(x,y,t)
&\geq& \frac{ V_{\mathbb{R}^n}(1)}{e(4\pi)^{n/2}}\left(\sqrt{\frac{d^2}{4t}+1}+\sqrt{\frac{d^2}{4t}} \right)^{-n}\exp\left[-\frac{d^2(x,y)}{4t}\right]\\
&&\times V_x^{-\frac{1}{2}}\left( \sqrt{t}\right)V_y^{-\frac{1}{2}}\left( \sqrt{t}\right)\nn
\end{eqnarray}
with $d=d(x,y)$. And from our Gaussian upper bound (\ref{GUBk=0}), we have
\begin{eqnarray}\label{HK-UB}
H(x,y,t)
\leq  2^{\frac{n}{2}} e^{2}G_{\max}V_x^{-\frac{1}{2}}\left(\sqrt{t}\right)V_y^{-\frac{1}{2}}\left(\sqrt{t}\right),
\end{eqnarray}
which follows from the function $G(x)=\left(\sqrt{1+x^2}+x\right)^ne^{-x^2}$ achieving its maximum value $G_{\max}=\left(\sqrt{n^2+1}+n\right)^{n/2}e^{-\frac{\sqrt{n^2+1}-1}{2}}$ at $x^2=\frac{\sqrt{n^2+1}-1}{2}$.\\

Given $0<\alpha<1$, for fixed $t>0$ and $y\in M$, set $u(x,s)=: H\left(x,y,s+\alpha t\right)$, which is a positive solution to the heat equation on $(x,s)\in M\times [0,\infty)$. Since we assume $Ric(M)\geq 0$, from Bishop's volume comparison Theorem \cite{BC}, 
\begin{equation}\label{VolumeComparison}
V_x\left(\sqrt{\alpha t+s}\right)\leq V_x\left(\sqrt{t}\right)\leq \alpha^{-\frac{n}{2}} V_x\left(\sqrt{\alpha t}\right),\; \forall\; 0\leq s\leq \left(1-\alpha\right)t.
\end{equation}

Define $A=2^{\frac{n}{2}} e^{2}G_{\max}V_x^{-\frac{1}{2}}\left(\sqrt{\alpha t}\right)V_y^{-\frac{1}{2}}\left(\sqrt{\alpha t}\right)$, the above inequality and the heat kernel upper bound (\ref{HK-UB}) imply 
$$u(x,s)\leq 2^{\frac{n}{2}} e^{2}G_{\max}V_x^{-\frac{1}{2}}\left(\sqrt{\alpha t+s}\right)V_y^{-\frac{1}{2}}\left(\sqrt{\alpha t+s}\right) \leq A,$$
for all $(x,s)\in M\times [0,\infty)$.\\

Thus, by Hamilton-Kotschwar's gradient estimate, Theorem \ref{H-gradient}, we have
$$s\left|\nabla \ln u(x,s)\right|^2\leq \ln\left(\frac{A}{u(x,s)}\right)=\ln\left(\frac{A}{H\left(x,y,s+\alpha t\right)}\right)$$
on $M\times \left[0,\left(1-\alpha\right)t\right]$. Evaluating at $s=\left(1-\alpha\right)t$ and applying our Gaussian lower bound (\ref{HK-LB}) and the volume comparison (\ref{VolumeComparison}), we conclude that
\begin{eqnarray}
&&\left(1-\alpha\right)t\left|\nabla \ln H\left(x,y,t\right)\right|^2\\
&\leq&\ln\frac{2^{n/2} e^{2}G_{\max}V_x^{-\frac{1}{2}}\left(\sqrt{\alpha t}\right)V_y^{-\frac{1}{2}}\left(\sqrt{\alpha t}\right)}{\frac{ V_{\mathbb{R}^n}(1)}{e(4\pi)^{n/2}}\left(\sqrt{\frac{d^2}{4t}+1}+\sqrt{\frac{d^2}{4t}} \right)^{-n}\exp\left[-\frac{d^2(x,y)}{4t}\right]V_x^{-\frac{1}{2}}\left( \sqrt{t}\right)V_y^{-\frac{1}{2}}\left(\sqrt{t}\right)}\nn\\
&\leq& 3+\frac{n}{2}\ln 8\pi +\ln \frac{G_{\max}}{V_{\mathbb{R}^n}(1)}-\frac{n}{2}\ln\alpha +n\ln\left(\sqrt{\frac{d^2}{4t}+1}+\sqrt{\frac{d^2}{4t}} \right) + \frac{d^2(x,y)}{4t} \nn\\
&=& C(n)-\frac{n}{2}\ln\alpha +n\ln\left(\sqrt{\frac{d^2}{4t}+1}+\sqrt{\frac{d^2}{4t}} \right) + \frac{d^2(x,y)}{4t}, \nn
\end{eqnarray}
with $C(n)=\frac{n}{2}\ln8\left(n+\sqrt{n^2+1}\right)+\ln \Gamma\left(\frac{n}{2}+1\right)+\frac{5-\sqrt{n^2+1}}{2} $, where we apply $V_{\R^n}(1)=\pi^{n/2}/\Gamma\left(\frac{n}{2}+1\right)$, where $\Gamma(\tau)$ is the Euler Gamma function.\\

And by Hamilton-Wu's Laplacian estimate, Theorem \ref{H-Laplacian}, we have
$$s\left(\frac{\Delta u(x,s)}{ u(x,s)}\right)\leq n+4\ln\left(\frac{A}{u(x,s)}\right)=n+4\ln\left(\frac{A}{H\left(x,y,s+\alpha t\right)}\right)$$
on $M\times \left[0,\left(1-\alpha\right)t\right]$. Evaluating at $s=\left(1-\alpha\right)t$ and applying our Gaussian lower bound (\ref{HK-LB}) and the volume comparison (\ref{VolumeComparison}), we conclude that
\begin{eqnarray}
&&\left(1-\alpha\right)t\left(\frac{\Delta H\left(x,y,t\right)}{H\left(x,y,t\right)}\right)
\leq n+4\ln\left(\frac{A}{H\left(x,y,t\right)}\right)\\\
&\leq& n+4C(n)-2n\ln\alpha+4n\ln\left(\sqrt{\frac{d^2}{4t}+1}+\sqrt{\frac{d^2}{4t}} \right) + \frac{d^2(x,y)}{t}. \nn
\end{eqnarray}
\qe\\

Next we prove Theorem \ref{Hamilton-HK-sharp} by choosing special $\delta$ in Theorem \ref{Hamilton-HK},

{\bf Proof of Theorem \ref{Hamilton-HK-sharp}:}
Set $\alpha(t,d)=\min\left\{\frac{1}{2},\left(\sqrt{1+\frac{d^2}{4t}}+\sqrt{\frac{d^2}{4t}}\right)^{-2}\right\}$ for any $t>0$ and $d\geq 0$, then we obtain
\begin{eqnarray*}
\left(1-\alpha(t,d)\right)&=&\max\left\{\frac{1}{2},\frac{2}{1+\sqrt{1+\frac{4t}{d^2}}}\right\},\\
-\ln\alpha(t,d) &\leq& \ln 2+2\ln\left(\sqrt{1+\frac{d^2}{4t}}+\sqrt{\frac{d^2}{4t}}\right).
\end{eqnarray*}
For any fixed $t>0$ and fixed two points $x, y \in M$, by taking $\delta=\delta(t,d_M(x,y))$ in Theorem \ref{Hamilton-HK}, both estimates of Theorem \ref{Hamilton-HK-sharp} follow.
\qe

\begin{Rem}
For the heat kernel of $\R^n$, $H(x,y,t)=\left(4\pi t\right)^{-n/2}e^{-\frac{d^2(x,y)}{4t}}$. As $\frac{d^2(x,y)}{4t}\rightarrow\infty$, the left-side of (\ref{HK-gradient-sharp}) is asymptotic to $\frac{d^2(x,y)}{4t}$, and the right-side  of (\ref{HK-gradient-sharp}) is asymptotic to $\frac{d^2(x,y)}{4t}+2n\ln\left(\sqrt{\frac{d^2}{4t}+1}+\sqrt{\frac{d^2}{4t}} \right)$,
which shows the sharpness of gradient estimate (\ref{HK-gradient-sharp}).\\
 
 And the Laplacian estimate (\ref{HK-Laplacian-sharp}) is sharp in the order of $\frac{d^2(x,y)}{4t}$ for  the heat kernel of $\R^n$, where the right-side of (\ref{HK-Laplacian-sharp}) is  asymptotically 4 times of the left-side  of  (\ref{HK-Laplacian-sharp}) as $\frac{d^2(x,y)}{4t}\rightarrow\infty$.\\
\end{Rem}

\section{\bf Manifolds with maximal volume growth}
In this section, as another application of our new lower bound estimate, Theorem \ref{Thm-GLB-delta}, we will give a simpler proof for the result concerning the asymptotic behavior in the time variable for the heat kernel as was proved in \cite{LiP-1} on a complete manifold $M$ with $Ric(M)\geq 0$ and maximal volume growth. If $B_x(r)$ denotes the geodesic ball of radius $r$ centered at $x\in M$, then we denote $V_x(r)$ and $A_x(r)$ to be the volume of $B_x(r)$ and the area of $\partial B_x(r)$, respectively. $M$ with maximal volume growth means that there exists $\theta > 0$ independent of $p\in M$, such that
\begin{equation}\label{MVG-rate}
\theta_p(r)=n^{-1}r^{1-n}A_x(r)\searrow \theta \quad \text{and}\quad r^{-n}V_p(r)\searrow \theta,\quad \text{as} \; r\rightarrow\infty.
\end{equation}
And it was proved in \cite{LiP} that for a fixed point $p \in M$, the function $t^{n/2}H(p, p, t)$ is monotonically non-decreasing with
$$t^{n/2}H(p, p, t)\nearrow \frac{V_{\R^n}(1)}{\left(4\pi\right)^{n/2}\theta}.$$

The following Lemma gives a slightly stronger result concerning the lower bound of the asymptotic
behavior in the time variable for the heat kernel $H(x,y,t)$ in \cite{LiP} and \cite{LTW},
\begin{Lemma}\label{liminf-asym-Lemma}
Let $M$ be a complete manifold with non-negative Ricci curvature and maximal volume growth. If $\gamma(t)=(y(t),t)$ is any path on $M\times(0,\infty)$ satisfying $d^2(x,y(t)) = o(t^2)$ as $t\rightarrow \infty$, then for any $x\in M$
\begin{equation}\label{liminf-asymptotic}
\liminf_{t\rightarrow \infty}V_x\left(\sqrt{t}\right)e^{\frac{d^2(x,y(t))}{4t}}H(x,y(t),t)\geq \frac{V_{\R^n}(1)}{\left(4\pi\right)^{n/2}}
\end{equation}
\end{Lemma}

\bp
From the lower bound estimate (\ref{GLB-delta}), for any $\delta>0$, we have
\begin{eqnarray*}
H(x,y(t),t)&\geq&e^{-\delta}\frac{ V_{\mathbb{R}^n}(R_{\delta}(t))}{V_x(R_{\delta}(t))}\left(4\pi t\right)^{-n/2}\exp\left[-\frac{d^2(x,y(t))}{4t}\right],
\end{eqnarray*}
with $d=d(x,y(t))$, and 
$$R_{\delta}(t)=\frac{\sqrt{d^2+4\delta t}-d}{2}=\frac{\delta\sqrt{ t}}{\sqrt{\frac{d^2}{4t}+\delta }+\sqrt{\frac{d^2}{4t}}}=\frac{2\delta}{1+\sqrt{1+\frac{4t}{d^2}\delta}} \frac{t}{d}.$$

{\bf Claim:} $\lim_{t\rightarrow\infty}R_{\delta}(t)=\infty$ under $d^2(x,y(t)) = o(t^2)$ as $t\rightarrow \infty$.\\

{\bf Case 1:} if $\limsup_{t\rightarrow\infty}\frac{d^2(x,y(t)}{4t}\leq A$ for some $0\leq A<\infty$,
\begin{eqnarray*}
\liminf_{t\rightarrow\infty}R_{\delta}(t)=\liminf_{t\rightarrow\infty}\frac{\delta\sqrt{ t}}{\sqrt{\frac{d^2}{4t}+\delta }+\sqrt{\frac{d^2}{4t}}}\geq \liminf_{t\rightarrow\infty}\frac{\delta\sqrt{ t}}{\sqrt{A+\delta }+\sqrt{A}}=\infty.
\end{eqnarray*}

{\bf Case 2:}  if $\limsup_{t\rightarrow\infty}\frac{d^2(x,y(t)}{4t}=\infty$, by passing sub-sequence
$t_i\rightarrow\infty$, without loss of generality, we may assume that $\lim_{t\rightarrow\infty}\frac{d^2(x,y(t)}{4t}=\infty$, then
\begin{eqnarray*}
\liminf_{t\rightarrow\infty}R_{\delta}(t)=\liminf_{t\rightarrow\infty}\frac{2\delta}{1+\sqrt{1+\frac{4t}{d^2}\delta}} \frac{t}{d}= \liminf_{t\rightarrow\infty}\delta \frac{t}{d}=\infty.
\end{eqnarray*}

Thus, using the maximal volume growth condition, we obtain
\begin{eqnarray*}
&&\liminf_{t\rightarrow \infty}V_x\left(\sqrt{t}\right)e^{\frac{d^2(x,y(t))}{4t}}H(x,y(t),t)\\
&\geq&\liminf_{t\rightarrow \infty}e^{-\delta}\frac{ V_{\mathbb{R}^n}(R_{\delta}(t))}{V_x(R_{\delta}(t))}\frac{V_x\left(\sqrt{t}\right)}{\left(4\pi t\right)^{n/2}}\\
&=&e^{-\delta}\frac{ V_{\mathbb{R}^n}(1)}{\theta}\frac{\theta}{\left(4\pi \right)^{n/2}}
=e^{-\delta}\frac{ V_{\mathbb{R}^n}(1)}{\left(4\pi \right)^{n/2}}
\end{eqnarray*}
Our estimate (\ref{liminf-asymptotic}) follows as $\delta\rightarrow 0$.

\qe

Next we give a simpler proof for the result concerning the asymptotic behavior in the time variable for the heat kernel as was proved in \cite{LiP},

\begin{Th}{\bf(Theorem 1 in \cite{LiP})}
Let $M$ be a complete manifold with non-negative Ricci curvature, $Ric(M)\geq 0$, and maximal volume growth. If $\gamma(t)=(y(t),t)$ is any path on $M\times(0,\infty)$ satisfying $d^2(x,y(t)) = o(t)$ as $t\rightarrow \infty$, then for any $x\in M$
\begin{equation}\label{lim-asymptotic}
\lim_{t\rightarrow \infty}V_x\left(\sqrt{t}\right)e^{\frac{d^2(x,y(t))}{4t}}H(x,y(t),t)= \frac{V_{\R^n}(1)}{\left(4\pi\right)^{n/2}}
\end{equation}
\end{Th}

\bp
For any $\delta>0$, using $R_{\delta}(t)=\frac{\sqrt{d^2+4\delta t}-d}{2}$ as defined in Theorem \ref{Thm-GLB-delta}, we'll estimate the upper bound of $V_x\left(\sqrt{t}\right)e^{\frac{d^2(x,y)}{4t}}H(x,y,t)$ as following:\\

When $0\leq d=d(x,y)\leq 2R_{\delta}(t)$, equivalent to $0\leq \frac{d^2}{4t}\leq \frac{\delta}{3}$. Applying the Harnack inequality (\ref{LYHarnack}) to $H(x,y,t)$ with $t_1=t$ and $t_2=2t$, 
\begin{eqnarray}\label{limsup-asym-est-0}
V_x\left(\sqrt{t}\right)e^{\frac{d^2(x,y)}{4t}}H(x,y,t)
& \leq & \frac{ V_x\left(\sqrt{t}\right)}{t^{n/2}}e^{\frac{d^2(x,y)}{2t}}(2t)^{n/2}H(x,x,2t)\\
&\leq& \frac{ V_x\left(\sqrt{t}\right)}{t^{n/2}} e^{2\delta/3}\frac{V_{\R^n}(1)}{\left(4\pi\right)^{n/2}\theta}.\nn
\end{eqnarray}

Hence combined with Lemma \ref{liminf-asym-Lemma}, since $\frac{ V_x\left(\sqrt{t}\right)}{t^{n/2}}\searrow \theta$ as $t\rightarrow\infty$, we obtain the asymptotic result (\ref{lim-asymptotic}) by send $\delta\rightarrow 0$, if  $d^2(x,y(t)) = o(t)$, as $t\rightarrow \infty$, as Theorem 1 in \cite{LiP}.

\qe

\end{document}